\documentclass[twocolumn]{autart}   
\let\oldlabel\label
\DeclareRobustCommand{\label}[1]{\oldlabel{b:#1}}
\usepackage{amsmath,amssymb,latexsym}
\usepackage{graphicx}       
\usepackage{enumitem}
\usepackage{here}
\usepackage{tikz}

\usepackage{cases}
\usepackage{url}

\let\oldref\ref
\renewcommand{\ref}[1]{\oldref{b:#1}}

\newcommand\norme[1]{\left \Vert #1 \right\Vert}

\newcommand\scal[2]{\left< #1,\, #2 \right>}

\newtheorem{theorem}{Theorem}
\newtheorem{lemma}[theorem]{Lemma}
\newtheorem{proposition}{Proposition}
\newtheorem{remark}{Remark}
\usepackage{framed}
\usepackage{comment}
\usepackage{xcolor}
\definecolor{shadecolor}{gray}{0.95}
\specialcomment{extended}{\begin{shaded}}{\end{shaded}}

\begin{document}

\begin{frontmatter}
\title{Lyapunov functions for linear damped wave equations in one-dimensional space with dynamic boundary conditions.} 

\author[Y]{Yacine Chitour}\ead{yacine.chitour@l2s.centralesupelec.fr},    
\author[H]{Hoai-Minh Nguyen}\ead{hoai-minh.nguyen@sorbonne-universite.fr},             
\author[C]{Christophe Roman}\ead{christophe.roman@lis-lab.fr} 

\address[Y]{Laboratoire des signaux et syst\`emes, Universit\'e Paris Saclay, Centralesupelec CNRS, Gif-sur-Yvette, France.}      
\address[H]{Laboratoire Jacques Louis Lions, 
Sorbonne Universit\'e, Paris, France.} 
\address[C]{Laboratoire informatique et syst\`eme, 
Aix-Marseille Universit\'e, Marseille, France.}

\begin{abstract}
{
This paper considers a one-dimensional wave equation on $[0,1]$, with dynamic boundary conditions of second order at $x=0$ and $x=1$, also referred to as Wentzell/Ventzel boundary conditions in the literature. 
In additions the wave is subjected to constant disturbance in the domain and at the boundary. 
This model is inspired by a real experiment.
By the means of a proportional integral control, the regulation with exponential converge rate is obtained when the damping coefficient is a nowhere-vanishing function of space. 
The analysis is based on the determination of appropriate Lyapunov functions and some further analysis on an associated error system. The latter  is proven to be exponentially stable towards an attractor.
Numerical simulations on the output regulation problem and additional results on related wave equations are also provided.
}

\end{abstract}
\end{frontmatter}

\maketitle
The wave equation is one of the classical partial differential equations.
The actual reason is that the wave equation is the continuous pendant of Newton's second law of motion, i.e., where 
momentum is equal to the sum of the forces.
As a consequence, it is also linked with the Euler-Lagrange framework, and therefore with the principle of least action.
For stationary systems, the energy is conserved and the action (or Lagrangian) is stationary.
Other physical phenomena are therefore associated with the wave equation such that electromagnetic law, and quantum phenomena with the Klein-Gordon equation.

In the control community, the wave equation has been mainly used for the modelization, estimation, and control of mechanical vibration and deformation phenomena. 
The regulation and control problem applied on the one-dimensional wave equation with dynamic boundary condition has attracted the attention of many researchers in the control community: 
crane regulation \cite{dandrea1992}, \cite{dandrea1994}, \cite{dandrea2000}, and \cite{conrad1998strong}, hanging cable immersed in water \cite{Bohn2013}, drilling torsional vibrations \cite{Saldivar2016}, \cite{terrandjeanne2019} ,\cite{Barreau2021}, \cite{Vanspranghe2022}, piezoelectric control \cite{Meurer2011}, and flexible structure \cite{Halevi2004}.
There are nowadays two main classes of issues : on the one hand, longitudinal variation with for example overhead crane and underwater cable, and, on the other hand, torsional variation with drilling string dynamics. 
The difference is on the control objective: one aims at controlling  the position in the first case, and instead the velocity in the second case.

The behavior of the wave equation is strongly related to its boundary conditions.
In the case of classical boundary condition (i.e., Dirichlet, Neumann, Robin) that issue is well understood  in the linear case and without high-order terms.  
Particular terms at one boundary can compensate for anti-damping terms at other boundaries and even in the domain, for example, see \cite{smyshlyaev2009boundary}, \cite{smyshlyaev2010boundary} and \cite{roman2018backstepping}.
Moreover, there are cases where even if the energy of the one-dimensional linear wave equation decreases along trajectories, it still does not decay exponentially \cite[Section 4]{li2017boundary}.

The wave equation under consideration is subject to two dynamic boundary conditions. This model results from an identification problem associated with a laboratory experiment \cite{roman2021}. 

\section{Problem statement.}
The considered system is defined for  $t\geqslant 0$ and for $x\in(0,1)$, by
\begin{subnumcases}{\label{sys_problem}}
v_{tt}(t,x)=(a(x)v_x)_x(t,x)-q (x)v_t(t,x)+f(x), \\
v_{tt}(t,1)=-\beta_1v_x(t,1)-\nu v_t(t,1)+U(t)+f_c,\quad \label{sys_problem_2}\\
v_{tt}(t,0) = \mu_1v_x(t,0)-\gamma_1 v_t(t,0)+f_{ac}, \label{sys_problem_3} \\
v(0,\cdot)=v_0,\quad v_t(0,\cdot)=v_1.
\end{subnumcases}
Here $U(t)$ is the control input and  we assume that
\begin{enumerate}[label=$(h_\emph{\arabic*})$]
\item\label{item1} the function $a:[0,1]\to \mathbb{R}_+^*$ is in $W^{1,\infty}(0,1)$ and that there exist 
$\underline{a}, \overline{a}>0$ such that $\underline{a}\leq a(\cdot) \leq \overline{a}$ a.e. on $[0,1]$. 
This function is associated with the mass and elasticity of the wave and it is also linked with the velocity.
\item\label{item2} The function $q:[0,1]\to \mathbb{R}_+^*$, describing the in-domain damping is in $L^{\infty}(0,1)$ and satisfies $\underline{q}\leq q(\cdot)\leq \overline{q}$ a.e. on $[0,1]$ for some $\underline{q}, \overline{q}>0$. 
\item\label{item3}  The constants $\beta_1$, $\gamma_1,\mu_1$ are positive real numbers, and $\nu $ is real.
\item\label{item4} The source terms $f(\cdot)$ is in $L^\infty (0,1)$, and the real constants $f_c,f_{ac}$ are unknown and therefore they cannot be used in the computation of the control law $U(t)$.
\end{enumerate}
The regularity of $a(\cdot)$ { stated in \ref{item1}} follows by classical arguments.
In detail, for the computation, we need $a(\cdot)v_x(t,\cdot)$ to be in $H^1(0,1)$. To be more precise everything will be the same as in the constant parameters case if $a(\cdot)v_x(t,\cdot)$ and $v_x(t,\cdot)$ have the same regularity. 
To get strong solutions for \eqref{sys_problem}, one needs to have that $ v\in H^2(0,1)$. Next, it can be easily shown that if $a\in W^{1,\infty}$ and $,  v\in H^2(0,1)$ then $a(\cdot)v_x(t,\cdot)$ is in $H^1(0,1)$.
Note that this is just a sufficient condition for the regularity.
We refer the reader to \cite[Chapter 21]{tartar2007introduction} for more details { about the regularity of $a$}.
In the sequel, we also need $qv_t(t,\cdot)^2$ need to be integrable, this means $q\in L^{\infty}(0,1)$. For $f$ we actually only need it to be integrable, it holds nonetheless $L^{\infty}(0,1)\subset L^1(0,1)$. 

{ The objective of the paper is to regulate $v_t(t,\cdot)$ to the constant reference value $v_1^{\text{ref}}$, by means of a proportional integral (PI) control law using the measurement of the velocity collocated with the actuation, $v_t(1,t)$, in other words, the control $U(t)$ can take the form
\begin{equation}
U(t):=-k(v_t(t,1)-v_1^{\text{ref}})-k_i\int_0^t (v_t(s,1)-v_1^{\text{ref}}) ds,\label{def_U}
\end{equation}
where the constants $k,k_i$ have to be chosen.
}
This can be equivalently written as
\begin{subnumcases}{\label{def}}
U(t)=-k(v_t(t,1)-v_1^{\text{ref}})-k_i\eta_v(t), \\ \dot{\eta}_v(t)= v_t(t,1)-v_1^{\text{ref}}, \quad \eta_v(0)=0. 
\end{subnumcases}
In the literature boundary conditions of the type \eqref{sys_problem_2}-\eqref{sys_problem_3} can be recast as Wentzell's boundary conditions \cite{fourrier2013regularity}. It involves a modification of the usual state space which in our case requires the addition of two finite-dimensional state variables, in a similar way as in \cite{Slemrod1989}, \cite{MIFDAL199737}, \cite{dandrea1994}, \cite{dandrea2000} and \cite{conrad1998strong}. When the wave equation is more than a one-dimensional, the reader is referred to \cite{fourrier2013regularity} { and \cite{BUFFE2017207}} and references therein.

{This type of control problem lies in robust output regulation. There has been an effort to extend the result and method from linear finite dimensional systems, to infinite dimensional systems.
We refer the reader to \cite{paunonen2010internal}, \cite{paunonen2015controller}, \cite{Paunonen2019}, more recently \cite{Vanspranghe2023} and reference within all of them.
These papers establish general results for example \cite{Vanspranghe2023} deals with non-linear systems. However, they are mostly based on either passivity, strong monoticity or exponential decay properties.
These properties often remain to be proven as it is the case of the present paper.
In \cite{terrand2019adding}, the authors establish general result on the PI control of infinite dimensional systems with the assumption beforehand on the exponential stability of the zero input system.}

{ The impact of the in-domain damping $q(x)$ can be an issue for the decay rate, as we can have some overdamping phenomenon. Intuitively, the damping should help the decay rate of the system. But as one can see in \cite{GUGAT201472} where a semi-linear wave equation is considered, the decay rate of the non-damped system is finite time, the addition of the damping degrade this performance to an exponential decay rate. Note that for the present case the in-domain damping is mandatory for the proof.}

{There are specific configurations of \eqref{sys_problem} that can be solved using more intricate and general control law designs, especially those tailored to address partial differential equations coupled with ordinary differential equations at the boundary. If disturbances are not considered, \cite{deutscher2018output} and its extension \cite{deutscher2019output} can be employed. Additionally, assuming that $a(\cdot)$ is constant allows the use of \cite{saba2017backstepping}, \cite{wang2020output}, or \cite{wang2020delay}.

These five papers primarily employ an infinite-dimensional backstepping approach, a development closely associated with the influential work of Miroslav Krstic \cite{krstic2008boundary}. Given that the wave equation can be expressed as a coupled heterodirectional hyperbolic partial differential equation (PDE), the main strategy in the aforementioned papers involves using backstepping transformations to decouple or cascade the PDE. This transforms the closed-loop system into a target system, the stability of which is easier to analyze. Notably, the uniqueness of the present paper lies in achieving exponential stability without the need for decoupling, thereby establishing new potential target systems for backstepping based design.
}

The closest approach associated with the present paper is \cite{terrandjeanne2019} where the velocity regulation with a PI is considered. However, the controlled boundary condition considered in \cite{terrandjeanne2019} is not a second order dynamic one, and thus is different from the one considered in this paper. Nevertheless, the boundary condition considered in \cite{terrandjeanne2019} implies the exponential stability even with small viscous anti-damping at the boundary opposite to the actuation. In the case under consideration, only viscous damping at the opposite boundary is considered, and exponential stability is achieved.
In \cite{conrad1998strong} the wave equation is subject to two dynamic boundary conditions.{ The authors establish asymptotic stability for the position stabilization and that the decay rate is not exponential}, and no viscous terms are considered for the zero input system. { In \cite{MIFDAL199737}, for the same model (as \cite{conrad1998strong}) the exponential stability toward the origin is obtained but the control law needs the knowledge of $v_{xt}(1,t)$. This can be related to the finite dimensional backstepping done in \cite{dandrea2000}.}
{ Studies have been conducted concerning the potential absence of exponential stabilization for wave-like equations, as evidenced by works such as \cite{morgul1994stabilization}, \cite{rao1995uniform}, and references therein. In a broader context, investigations into this issue extend to more general setups, as seen in \cite{gibson1980note}, \cite{triggiani1989lack}, and related references.}

PI controllers have been successfully and recently used in order to regulate linear and non-linear PDE, see \cite{coron2019pi}, \cite{lhachemi2022proportional}. An identification procedure has been presented in \cite{roman2021} for the system \eqref{sys_problem} without source terms on experimental data. This means that the considered problem can be associated with an experimental setup. A first study has been made on this system in \cite{roman:hal-03902877} using classical form a Lyapunov functional but it failed to prove the exponential stability. Only asymptotic stability was established, by using the LaSalle invariance principle.

This paper provides a new term in the Lyapunov functional and an associated methodology, for the present setup. The proof of the exponential stability is given in Section~\ref{sec:proof:main}. In Section~\ref{sec:next_result}, this proof is compared with existing results. Next the proof of the robustness of the controlled system is given in Section~\ref{sec:proof:robust}. { Then in Section~\ref{sec:zero_input} we study, using the same approach, simpler cases where one boundary condition is a Dirichlet one and this allows us to establish the exponential stability of the zero input system in the undisturbed case.} The last part of the paper deals with numerical simulations. The numerical scheme is not derived from the usual approximation of space and time derivatives. We used the fact that the wave equation can be derived from the Lagrangian and the least action principle to approximate the system space energy by a finite dimensional continuous time Euler-Lagrange equation. The finite dimensional continuous time system is then numerically solved by using symplectic integrators. This suggested numerical scheme is new up to the authors' knowledge and provides an interesting alternative compares to more standard discretization schemes.

{\bf Notations:} 
If $I$ is an interval of real numbers, $L^2(I;\mathbb{R})$ denotes (the class of equivalence of) square-integrable functions from $I$ to $\mathbb{R}$. Moreover $L^2([0,1];\mathbb{R})$ is abusively denoted $L^2(0,1)$. Furthermore $H^n$ denotes the Sobolev space $W^{n,2}$, i.e., 
\begin{equation}
u\in H^1 \Leftrightarrow u\in L^2,\ u'\in L^2,
\end{equation}
in which $u'$ denotes the derivative of $u$.

\section{Main result}
{
To achieve our objective, we perform a change of variable in order to obtain an error variable $u(\cdot,\cdot)$ and to prove exponential decay of its partial derivatives. 

The error variable $u(\cdot,\cdot)$ is defined as follows,}
for every $(x,t)\in [0,1]\times [0,\infty)$
\begin{align}
u(t,x)&:= v(t,x)-tv_1^{\text{ref}} \notag\\&\quad +\int_0^x\frac{1}{a(s)}\int_0^s [-v_1^{\text{ref}}q(\chi)+f(\chi)]d\chi ds \notag\\ +&\frac{a(0)}{\mu_1}[-\gamma_1 v_1^{\text{ref}}+f_{ac}]\int_0^x\frac{1}{a(s)}ds ,\label{def_u} \\
\eta_2(t)&:=\eta_v(t) - \frac{\beta_1 }{{k_i} a(1)}\int_0^1[-v_1^{\text{ref}}q(s)+f(s)]ds\notag\\&\quad -\frac{\beta_1 a(0)}{{k_i} \mu_1 a(1)}[-\gamma_1 v_1^{\text{ref}}+f_{ac}]+\frac{\nu  v_1^{\text{ref}}-f_c}{{k_i}}.
\end{align}
{
Note that, for every $(x,t)\in [0,1]\times [0,\infty)$,
\begin{align}
 u(t,x)&= v(t,x)-tv_1^{\text{ref}}+F(x),\\
u_t(t,x)&=v_t(t,x)-v_1^{\text{ref}},\label{eq:u-v}
\end{align}
where we have gathered all the uncertainties in the function $F(\cdot)$ and it is immediate to deduce from \eqref{eq:u-v} that proving exponential decay of $u_t$ (in an appropriate sense) is equivalent to prove it for $v_t-v_1^{\text{ref}}$ and hence to achieve the desired control objective.

From now, we will therefore focus on the error variable $u(\cdot,\cdot)$.
}Direct computations yield that it is the solution of the following system:
\begin{subnumcases}{\label{sys:propW2-W1}}
u_{tt} - (a(x)u_x)_x  = -q(x)u_t,\notag\\ \quad (t,x)\in \mathbb{R}^+ \times (0, 1),  \\
u_t(t,1)=\eta_1(t)\\
u_t(t,0)=\xi_1(t)\\
\dot\eta_1(t)=-\alpha_1\eta_1(t)-\alpha_2\eta_2(t)-\beta_1u_x(t,1),\label{eq:xi1}\\
\dot\eta_2(t)=\eta_1(t),\label{eq:eta}\\
\dot\xi_1(t)=-\gamma_1\xi_1(t)+\mu_1u_x(t,0),\label{eq:xi0}\\
u(0, \cdot) = u_0, \quad u_t (0, \cdot) = u_1 \quad \mbox{ on } (0, 1), \\
\eta_1(0)=\eta_0,\quad \eta_2(0)=\eta_{2,0} \quad \xi_1(0)=\xi_0.
\end{subnumcases}
where {$\alpha_1:=k+\nu $ and $\alpha_2:=k_i$}, and $k$ is chosen such that $\alpha_1$ is positive. 

Consider the following Hilbert spaces
\begin{align}
\quad X_w:&=H^1((0,1);\mathbb{R})\times L^2((0,1);\mathbb{R})\times \mathbb{R}^3,\label{def:Xw}\\
X_s:&=H^2((0,1);\mathbb{R})\times H^1((0,1);\mathbb{R})\times \mathbb{R}^3. \label{def:SS}
\end{align}
The wave equation is associated with the following abstract problem
\begin{subnumcases}{\label{sys:abs}}
\dot{\mathcal{X}}(t)+\mathcal{A}\mathcal{X}(t)=0,\\
\mathcal{X}(0)=\mathcal{X}_0\in \text{Dom}(\mathcal{A})\subset X_s \subset X_w,
\end{subnumcases}
in which 
\begin{align}
\forall z \in \text{Dom}(\mathcal{A}),\ \mathcal{A}z:=\begin{bmatrix}
  -z_2\\-(az_1')'+qz_2 \\ \alpha_1z_3+\alpha_2z_4+\beta_1z_1'(1) \\ -z_3\\\gamma_1z_5-\mu_1z_1'(0) 
\end{bmatrix},\label{def_A}
\end{align}
and 
\begin{align}
\text{Dom}(\mathcal{A}):=\{ z\in X_s;\  z_2(1)=z_3,\, z_2(0)=z_5\}.
\end{align}
Our well-posed result goes as follows.

{
\begin{theorem}\label{th:wellposeness} 
Considering assumption \ref{item1} and \ref{item2}, the abstract problem \eqref{sys:abs} is well-posed. In order words for any initial data $\mathcal{X}_0\in \text{Dom}(\mathcal{A})$,
there exists a unique solution to the abstract problem \eqref{sys:abs}, such that for any $t\geq 0$, $\mathcal{X}(t)\in \text{Dom}(\mathcal{A})\subset X_s$ and 
\begin{align}
  \mathcal{X} \in C^0([0,\infty);\text{Dom}(A)) \cap C^1([0,\infty);X_w),
\end{align}
$X_w$ is the state space of weak solutions and the Hilbert space considered and is defined in \eqref{def:Xw}. $X_s$ is the state space of strong solutions and is defined in \eqref{def:SS}.

In addition, for all initial data $\mathcal{X}_0 \in X_w$, there exists a weak solution
$\mathcal{X}(t) \in X_w$ to the abstract problem~\eqref{sys:abs} given by
\begin{align}
\mathcal{X}(t) = S(t)\mathcal{X}_0,  
\end{align}
in which $S$ is the $C_0$-semigroup generated by the unbounded operator $A$. Moreover, it holds
\begin{align}
  \mathcal{X} \in C^0([0,\infty);X_w).
\end{align}
\end{theorem}}

The proof is based on finding a transformation such that the abstract problem is associated with a linear maximal monotone operator.
Then the conclusion is drawn by using the Hille-Yosida theorem.
{ The part on weak solutions holds true from the fact that $\text{Dom}(\mathcal{A})$ is dense in $X_w$, and therefore $S(t)$ defined a strongly continuous map from $X_w$ to $X_w$}. Details are provided in Appendix~\ref{app_well_posed}.
The state is 
\begin{align}
  \mathcal{X}(t) := [&u(t,\cdot),\, u_t(t,\cdot),\notag \\ &\, \eta_1(t),\, \eta_2(t),\, \xi_1(t)] \in \text{Dom}(\mathcal{A}) \subset X_s.
\end{align}

We define the \emph{energy} $E_u$ of a solution of \eqref{sys:propW2-W1} as $\forall t\geq 0$
\begin{align}
E_u(t):=\frac{1}{2}\int_0^1(u_t(t,x)^2+a(x)u_x(t,x)^2)\, dx.
\end{align}

Note that this energy is invariant by translations with constants, i.e., $E_u=E_v$ if $u-v$ is a constant function. 
Moreover, the absolutely continuous function \mbox{$u(\cdot,1)-\eta_2(\cdot)$} is constant along a trajectory of \eqref{sys:propW2-W1} and equal to $u_*$ where 
\begin{align}
u_*:=u_0(1)-\eta_2(0).  \label{eq:u*}
\end{align}

Our objective is to establish the exponential stability of the trajectory with respect to the following attractor 
\begin{align}
    S:=\{&z\in X_w,\ z_1(\cdot)=d,d\in\mathbb{R},\, z_2(\cdot)=0,\notag\\ & z_3 = 0,\, z_4=0,\, z_5=0\}. 
\end{align}
This attractor is the kernel of the following functional 
\begin{align}
\Gamma(\mathcal{X}(t)):=& \int_0^1 [u_t^2(t,x)+u_x^2(t,x)] dx \notag\\& + \eta_1^2(t)+\eta_2^2(t)+\xi_1^2(t), \label{def_Gamma}
\end{align}
indeed it holds
\begin{align}
\Gamma(z)=0 \Leftrightarrow z\in S.
\end{align}
We establish the following result.
\begin{theorem}\label{th:main1}
Consider the 1D wave equation \eqref{sys:propW2-W1} with the assumptions \ref{item1}, \ref{item2}, \ref{item3}, and with $\alpha_2,\alpha_1>0$.  Then, there exist a positive constant $\rho$, and a positive constant $M$ 
such that, for every weak solution $\mathcal{X}$, it holds, 
\begin{align}
    \Gamma(\mathcal{X}(t))\leqslant M\Gamma(\mathcal{X}(0))e^{-\rho t}.\label{Gamma_theo}
\end{align}
and the system is exponentially stable towards the attractor $S$.

In addition it holds that $\max_{x\in[0,1]}\vert u(t,x)-u_*\vert$ tends exponentially to zero as $t$ tends to infinity, with a decay rate larger than or equal to $\rho$.

\end{theorem}

{
\begin{theorem}\label{theo_robust}
Under assumption \ref{item1}-\ref{item2}, and for any $k_i=\alpha_2>0$ and $\beta_1,\mu_1>0$, the conclusion on 
Theorem~\ref{th:main1} still holds if
\begin{align}
\alpha_1 > \frac{-\beta_1\kappa_2}{2a(1)}+\kappa_1\kappa_2,\\
\gamma_1 > \frac{-\mu_1\kappa_2}{2a(0)}+\kappa_1\kappa_2,
\end{align}
where 
\begin{align}
\kappa_1:= \frac{1}{\underline{a}}\big( 3\overline{a}+\frac{\Vert a_x\Vert_{L^\infty}}{2} + \frac{\overline{q}}{2} \big), \label{def_kappa_1}\\
\kappa_2:= \frac{2}{\underline{q}}(1+\frac{\overline{q}}{2}+\kappa_1).\label{def_kappa_2}
\end{align}

\end{theorem}
\begin{remark}
    The link between $\alpha_1$ and $\nu$ is defined right below \eqref{sys:propW2-W1}.
This theorem means in particular that \eqref{Gamma_theo} 
holds in a robust way and the regulation can even admit small anti-damping at the uncontrolled boundary for certain value of $\mu_1$, $a(\cdot)$, and $q(\cdot)$.
\end{remark}

}

  

\section{Proof of Theorem~\ref{th:main1}\label{sec:proof:main}}
This proof follows a standard strategy: the result is first established for strong solutions by the determination of Lyapunov functions verifying an appropriate differential inequality, and then it is extended to weak solutions by a classical density argument. Hence, in the sequel, solutions of \eqref{sys:propW2-W1} are all assumed to be strong.

We start with the time derivative of $E_u$ along a strong solution. It holds for $t\geq 0$
\begin{align}
\label{eq:dot-E0}
\dot E_u&=-\int_0^1qu_t^2\, dx+a(1)\eta_1(t)u_x(t,1)\notag\\&-a(0)\xi_1(t)u_x(t,0).  
\end{align}
One also has, for $t\geq 0$, after using \eqref{eq:xi1} and \eqref{eq:eta}
\begin{align}
a(1)\eta_1(t)u_x(t,1)&=-\frac{a(1)}{\beta_1}\eta_1(t)\Big(\dot\eta_1(t)+\alpha_1\eta_1(t)\notag\\&+\alpha_2\eta_2(t)\Big)\nonumber\\
&=-\frac{d}{dt}\Big(\frac{a(1)}{2\beta_1}(\eta_1^2(t)+\alpha_2\eta_2^2(t))\Big)\notag\\&-\frac{a(1)\alpha_1}{\beta_1}\eta_1^2(t).\label{eq:xi1ux}
\end{align}
Similarly, one also has, for $t\geq 0$, after using \eqref{eq:xi0}
\begin{align}
-a(0)\xi_1(t)&u_x(t,0)=-\frac{a(0)}{\mu_1}\xi_1(t)\Big(\dot\xi_1(t)+\gamma_1\xi_1(t)\Big)\nonumber\\
&=-\frac{d}{dt}\Big(\frac{a(0)}{2\mu_1}\xi_1^2(t)\Big)-\frac{a(0)\gamma_1}{\mu_1}\xi_1^2(t).\label{eq:xi0ux}
\end{align}
Define for $t\geq 0$
\begin{align}
\label{eq:E1}
F(\mathcal{X}(t)):=E_u(t)+\frac{a(1)}{2\beta_1}\eta_1^2(t)+\frac{a(0)}{2\mu_1}\xi_1^2(t).  
\end{align}
Then, by gathering \eqref{eq:dot-E0}, \eqref{eq:xi1ux} and \eqref{eq:xi0ux}, one deduces that, for $t\geq 0$,
\begin{align}
\frac{d}{dt}\Big(F&+\frac{a(1)\alpha_2}{2\beta_1}\eta_2^2(t)\Big)=
-\int_0^1qu_t^2\, dx\notag\\& -\frac{a(1)\alpha_1}{\beta_1}\eta_1^2(t)-\frac{a(0)\gamma_1}{\mu_1}\xi_1^2(t).\label{eq:E_00}
\end{align}
 
{ To conclude on the exponential stability we also need a negative term in $u_x^2$ and $\eta_2^2$.}
We next consider an extra term which will be added in the candidate Lyapunov function in the sequel. From \eqref{eq:u*} it holds that 
\begin{align}
\eta_2(t)=u(t,1)-u_*,\quad t\geq 0.  
\end{align}
Set
\begin{align}
\xi_2(t):=u(t,0)-u_*,\quad t\geq 0.  
\end{align}
One has, for $t\geq 0$, that 
\begin{align}
\frac{d}{dt}&\Big(\int_0^1(u-u_*)u_t\, dx\Big)=\int_0^1u_t^2+\int_0^1(u-u_*)u_{tt}\nonumber\\
=&\int_0^1u_t^2+\int_0^1(u-u_*)(au_x)_x-\int_0^1q(u-u_*)u_t\nonumber\\
=&\int_0^1u_t^2-\int_0^1au_x^2-\frac{d}{dt}\Big(\int_0^1\frac{q}2(u-u_*)^2\, dx\Big)\nonumber\\
&+a(1)\eta_2u_x(t,1)-a(0)\xi_2u_x(t,0).\label{eq:uut}
\end{align}
Using \eqref{eq:xi1} and \eqref{eq:xi0}, one deduces after computations similar to those performed to get \eqref{eq:xi1ux} and \eqref{eq:xi0ux}, that, for $t\geq 0$,
\begin{align}
\eta_2u_x(t,1)=&\frac{-\alpha_2\eta_2^2(t)+\eta_1^2(t)}{\beta_1}\notag\\&-\frac{d}{dt}\Big(\frac{\frac{\alpha_1}2\eta_2^2(t)+\eta_1(t)\eta_2(t)}{\beta_1}\Big),\label{eq:eta1xi1}\\
-\xi_2u_x(t,0)=&\frac{\xi_1^2(t)}{\mu_1}-\frac{d}{dt}\Big(\frac{\frac{\gamma_1}2\xi_2^2(t)+\xi_2(t)\xi_1(t)}{\mu_1}\Big)\label{eq:eta1xi0}.
\end{align}
We next define for $t\geq 0$
\begin{align}
W(\mathcal{X}(t))&=\int_0^1(u-u_*)u_t\, dx+\int_0^1\frac{q}2(u-u_*)^2\, dx
  \notag\\ +&\frac{a(1)}{\beta_1}\Big(\frac{\alpha_1}2\eta_2^2(t)+\eta_2(t)\eta_1(t)\Big)\nonumber\\
  +&\frac{a(0)}{\mu_1}\Big(\frac{\gamma_1}2\xi_2^2(t)+\xi_2(t)\xi_1(t)\Big).\label{eq:W0}
\end{align}  
Gathering \eqref{eq:uut}, \eqref{eq:eta1xi1} and \eqref{eq:eta1xi0}, it holds for $t\geq 0$
\begin{align}
\dot W=&-\int_0^1au_x^2-\frac{a(1)\alpha_2}{\beta_1}\eta_2^2(t)+\int_0^1u_t^2
  \notag\\+&\frac{a(1)}{\beta_1}\eta_1^2(t)+\frac{a(0)}{\mu_1}\xi_1^2(t).
  \label{eq:deuxieme}
\end{align}
We finally define the candidate Lyapunov function $V$ used for proving Theorem~\ref{th:main1}, which is positive definite for some constant $\ell$ such that $\sqrt{2\underline{q}}>\ell>0$ by
\begin{align}\label{eq:vu}
V(\mathcal{X}(t))=&F(\mathcal{X}(t))+\frac{a(1)\alpha_2}{2\beta_1}\eta_2^2(t)\notag\\&+\ell W(\mathcal{X}(t)),\quad \geq 0.
\end{align}
Putting together \eqref{eq:E1} and \eqref{eq:W0}, it holds for $t\geq 0$,
\begin{align}  
V(\mathcal{X}(t))=&E_u(t)\notag\\&+\ell\int_0^1\Big((u-u_*)u_t+\frac{q}2(u-u_*)^2\Big)\, dx\nonumber\\
&+\frac{a(1)}{2\beta_1}\Big(\eta_1^2+{\alpha_2\eta_2^2+\ell(2\eta_2\eta_1+\alpha_1\eta_2^2)}\Big)
\nonumber\\&+
\frac{a(0)}{2\mu_1}\Big(\xi_1^2+\ell(2\xi_2\xi_1+\gamma_1\xi_2^2)\Big).\label{eq:vu2}
\end{align}
and similarly, putting together \eqref{eq:E_00} and \eqref{eq:deuxieme},
it holds for $t\geq 0$,
\begin{align}  
\dot V=&-\int_0^1(\frac{q}2-\ell)u_t^2\, dx-\ell\int_0^1au_x^2\, dx\nonumber\\
&-\frac{a(1)}{\beta_1}\Big((\alpha_1-\ell)\eta_1^2+\alpha_2\ell\eta_2^2\Big)\notag\\&-\frac{a(0)(\gamma_1-\ell)}{\mu_1}\xi_1^2.\label{eq:dotvu2}
\end{align}
{ The purpose of $V$ defined in \eqref{eq:vu2} compared with $F$ is to make negative terms in $u_x^2$ and $\eta_2^2$ appear. Next we compare the functional $V$ to the functional $\Gamma$ defined in \eqref{def_Gamma}.}

\begin{proposition}\label{prop1}
With the notations above, and $\Gamma$ defined in \eqref{def_Gamma}, there exist $\ell>0$ and two positive constants $c,C,\rho>0$ such that for every strong solution $\mathcal{X}(t)$ of \eqref{sys:propW2-W1}, one gets, for $t\geq 0$, 
\begin{align}
c \Gamma (\mathcal{X}(t))&\leq V(\mathcal{X}(t))\leq C\Gamma(\mathcal{X}(t)),\label{eq:ET-est0}\\
\dot V (\mathcal{X}(t))&\leq-C\rho \Gamma(\mathcal{X}(t)).\label{eq:lyap}
\end{align} 
\end{proposition}
{
\begin{remark}
  Using $\alpha_1$ and $\alpha_2$ as tuning parameters one can show that a necessary condition for 
  \begin{align}
    \dot V(\mathcal{X}(t)) \leqslant - C\rho \Gamma (\mathcal{X}(t)) .
  \end{align}
  is that 
  \begin{align}
    C\rho < \min \{ \frac{\underline{a}\underline{q}}{4},\ \underline{a}\sqrt{2\underline{q}},\ \frac{\underline{a} a(0)\gamma_1}{a\mu_1+a(0)} \}.
  \end{align}
  This upper bound is deduced from the next inequalities extracted from \eqref{eq:dotvu2} and the condition for $V$ to be definite positive. 
  \begin{align}
    \sqrt{2\underline{q}}>\ell \\
    \frac{q}{2}-\ell > C\rho\\
    \ell a >C\rho \\
    \frac{a(0)(\gamma_1-\ell)}{\mu_1}> C\rho
  \end{align}
  Moreover as $c$ and $C$ does not depend on $q$. It holds for the decay rate $\rho$
  \begin{align}
    \rho \underset{\underline{q}\to  0}\longrightarrow 0.
  \end{align}
  The suggested approach allows us only to conclude for stability when $q=0$, and in this case we can stop at \eqref{eq:E_00}. Nevertheless following \cite{roman:hal-03902877} or \cite{conrad1998strong} we could use LaSalle's invariance principle to establish asymptotic stability. If in addition $\alpha_2=0$, in the case of no integrator the system falls as a one-dimensional particular case of \cite[Theorem 1.2]{BUFFE2017207}, and therefore the decay rate is at least logarithmic.
\end{remark}}
\begin{pf}
  { Using \eqref{eq:u*} and \eqref{eq:E1} one can observe} that for every $t\geq 0$ and $x\in [0,1]$ it holds
\begin{align}
\label{eq:final}
\vert u(t,x)-u_*\vert^2&\leq 2\vert u(t,x)-u(t,1)\vert^2+2\eta_2^2(t)\notag\\&\leq 
2\int_0^1u_x^2(t,x)\, dx+2\eta_2^2(t)\notag\\&\leq \frac{4}{\underline{a}}E_u(t)+2\eta_2^2(t),  
\end{align}
As an immediate consequence, one gets that, for $t\geq 0$,
\begin{align}
\int_0^1(u-u_*)^2\, dx&\leq \frac{4}{\underline{a}}E_u(t)+2\eta_2^2(t),\label{eq:finalu}\\
\xi_2^2(t)&\leq \frac{4}{\underline{a}}E_u(t)+2\eta_2^2(t).\label{eq:finalxi2}
\end{align}
The proof of \eqref{eq:ET-est0} relies now on the combination of \eqref{eq:vu2}, \eqref{eq:finalu} and \eqref{eq:finalxi2}, several completions of squares and the Cauchy-Schwartz inequality. As for the argument of \eqref{eq:lyap}, it is obtained similarly by using \eqref{eq:dotvu2}, \eqref{eq:finalu} and \eqref{eq:finalxi2},
\end{pf}

Relying on Proposition~\ref{prop1}, we complete the proof of Theorem~\ref{th:main1}.

From \eqref{eq:ET-est0} and \eqref{eq:lyap}, it follows that $\dot V\leq -\rho V$ hence yielding exponential decrease of $V$ at the rate $\rho$ and the similar conclusion holds for $\Gamma$, thanks to \eqref{eq:lyap}. All items of Theorem~\ref{th:main1} are proven 
after using \eqref{eq:final}.

\section{Discussion on the proof of the Theorem~\ref{th:main1} \label{sec:next_result}}

There exist cases where the linear one-dimensional wave does not decay exponentially. For example, the solution $u$ of the system
\begin{subnumcases}{\label{sys_1}}
u_{tt}(t,x)=u_{xx}(t,x), \quad \in \mathbb{R}^+\times (0,1),\\
u(t,0)=0,\\
u_{tt}(t,1)=-u_x(t,1)-u_t(t,1),
\end{subnumcases}
does not decrease exponentially towards the origin, see \cite[Section 4]{li2017boundary}. { It follows a $t^{-1}$ sharp decay rate. The addition/suppression of one term can make the decay rate drastically different, for example
\begin{subnumcases}{\label{sys_1a}}
  u_{tt}(t,x)=(au_{x})_x(t,x), \quad \in \mathbb{R}^+\times (0,1),\\
  u_{tt}(t,0)=u_x(t,0),\\
  u_{tt}(t,1)=-u_x(t,1)-u_t(t,1)-u(t,1)\notag\\ {\color{white}u_{tt}(t,1)=} -u_{xt}(t,1),
  \end{subnumcases}
is exponentially stable \cite{MIFDAL199737}, whereas 
\begin{subnumcases}{\label{sys_1b}}
u_{tt}(t,x)=(au_{x})_x(t,x), \quad \in \mathbb{R}^+\times (0,1),\\
u_{tt}(t,0)=u_x(t,0),\\
u_{tt}(t,1)=-u_x(t,1)-u_t(t,1)-u(t,1),
\end{subnumcases}
is not exponentially stable, see \cite{conrad1998strong}. However the solution of \eqref{sys_1a} need to be more regular, see \cite{MIFDAL199737}.} The energy of the following two systems
\begin{subnumcases}{}
u_{tt}(t,x)=u_{xx}(t,x), \quad \in \mathbb{R}^+\times (0,1),\\
u_x(t,0)=u_t(t,0),\label{sys_2:2}\\
u_x(t,1)=-u_t(t,1), \label{sys_2:3}
\end{subnumcases}
and
\begin{subnumcases}{}
u_{tt}(t,x)=u_{xx}(t,x), \quad \in \mathbb{R}^+\times (0,1),\\
u_x(t,0)=u_t(t,0)\label{sys_3:2},\\
u_{tt}(t,1)=-u_x(t,1)-u_t(t,1),\label{sys_3:3}
\end{subnumcases}
are exponentially decreasing \cite{roman2016}. Typically, for both previous cases, the exponential decrease and stability can be obtained via Energy/Lyapunov approach using cross terms in the following form.
\begin{align}
\int_{0}^{1}(1+x)u_tu_x dx,\label{u_tu_x}
\end{align}
which can make negative term as $u_x^2$ and $u_t^2$ appear for the Energy/Lyapunov functional derivative. This pervious term implies boundary terms in the following form
\begin{align}
\left[u_x^2+u_t^2\right]_{x=0}^1, \label{buxut}
\end{align}
{in the case of \eqref{sys_2:2}-\eqref{sys_2:3} or \eqref{sys_3:2}-\eqref{sys_3:3} we { can manage to handle}  this term. However, this is problematic when considering both boundary conditions as \eqref{eq:xi1} and \eqref{eq:xi0}. Indeed, even when $\alpha_2=0$, we do not arrive to cope with the term $u_x^2$ both in $0$ and $1$. This incapacity to handle the term $u_x^2$ with both dynamic boundary conditions is properly shown in \cite{roman:hal-03902877} with more general form of $u_xu_t$ cross terms, and considering a large family of reformation as hyperbolic PDE for example.\par

In particular the term \eqref{buxut} can be also taken care of if we have damped position terms ($u$) on the domain. Indeed in this case, this term enable us to use cross terms like 
\begin{align}
\int_0^1 uu_t,\label{buut}
\end{align}
The exponential stability of the linear wave equation at the origin with both dynamic boundary condition and damped in velocity and position everywhere is established in \cite[Chapter 9]{roman2018boundary}. We stress that the paper deals with velocity regulation which has been transformed to velocity exponential stability.}
The term \eqref{buut} is close to the one we suggest 
\begin{align}
\int_0^1 (u-u_*)u_t, 
\end{align}
This mostly corresponds to the beforehand knowledge of the limit value of $u$ for the system. This can be made because the integrator part of the system captures the distance between the state and the attractor. In our case this term can be added because $q$ is strictly positive, see \eqref{eq:uut}.

{
\section{Proof on Theorem~\ref{theo_robust}\label{sec:proof:robust}}
We start from the proof of Theorem~\ref{th:main1}, in \eqref{eq:dotvu2}, then we compute the derivative of the following cross term, using integration by parts
\begin{align}
  \frac{d}{dt}&\Big( \int_0^1(1-2x)u_xu_{t}\, dx\Big)=-\frac{a(1)u_x^2(t,1)+a(0)u_x^2(t,0)}2\notag\\&
  +\int_0^1(2a+\frac{(a(1-2x))'}2)u_x^2\, 
  -\frac{\eta_1^2(t)+\xi_1^2(t)}2dx\notag\\&
  +\int_0^1u_t^2\, dx 
  -\int_0^1(1-2x)qu_xu_{t}\, dx.\label{eq:uxut}
\end{align}

The above cross term can be used to make negative terms in $\eta_1^2$ and $\xi_1^2$ appear at the cost of positive terms in $u_t^2$ and $u_x^2$.

Consider that $\vert \ell_2 \vert <\sqrt{\underline{a}}$ then
\begin{align}
 G_u=V_u+\ell_2 \int_0^1(1-2x)u_xu_{t}\, dx,
\end{align}
is positive. 

Gathering \eqref{eq:dotvu2} and \eqref{eq:uxut}, and using the Young's inequality, the derivative of $G_u$ along the trajectory is
\begin{align}  
 \dot G_u\leq &-\int_0^1(\frac{q}2-\ell-\ell_2-\ell_2\frac{q}{2})u_t^2\, dx\notag\\& -\int_0^1(\ell a -\ell_2(2a+\frac{(a(1-2x))'}2)-\ell_2\frac{q}{2})u_x^2\, dx\nonumber\\
 -& \Big(\frac{a(1)}{\beta_1}(\alpha_1-\ell)+\frac{\ell_2}{2}\Big)\eta_1^2-(\frac{a(0)(\gamma_1-\ell)}{\mu_1}+\frac{\ell_2}{2})\xi_1^2\, 
\notag\\& -\frac{a(1)}{\beta_1}\alpha_2\ell\eta_2^2
\end{align}
The exponential stability still holds if the following inequalities hold
\begin{align}
\frac{\underline{q}}2-\ell-\ell_2(1+\frac{\overline{q}}{2})>0,\label{eq:ineq1}\\  
 \ell \underline{a} -\ell_2(3\overline{a}+\frac{\overline{a'}}{2}+\frac{\overline{q}}{2}) >0,\label{eq:ineq2} \\ 
 2a(1)(\alpha_1-\ell)+\beta_1 \ell_2>0, \label{eq:ineq3}\\
2a(0)(\gamma_1-\ell)+\mu_1 \ell_2 >0. \label{eq:ineq4}
\end{align}
A sufficient condition for the four previous inequalities to hold is
\begin{align}
\ell_2<\kappa_2,\\  
\ell >\kappa_1\kappa_2, \\ 
2a(1)(\alpha_1-\kappa_1\kappa_2)+\beta_1\kappa_2>0, \\
2a(0)(\gamma_1-\kappa_1\kappa_2)+\mu_1\kappa_2 >0. 
\end{align} 
where $\kappa_1$ and $\kappa_2$ are defined in \eqref{def_kappa_1}-\eqref{def_kappa_2}. This concludes the proof.
}

\section{Exponential stability for the zero input system with no disturbance.\label{sec:zero_input} }

In the following we investigate and establish results on associated problems. 
{ We first start with a wave equation subject to a Dirichlet's boundary conditions and a 2nd order dynamic boundary condition. The second system we add an integral action to the dynamics boundary condition. The third and last system consist of a wave equation with both 2nd order dynamics boundary conditions, and correspond to the zero input system with no disturbance.}

\begin{proposition}\label{prop:W1-D}
Consider the following 1D wave equation  
\begin{subnumcases}{\label{sys_prop_2}}
u_{tt} - (au_x)_x  = -qu_t, \quad (t,x)\in \mathbb{R}^+ \times (0, 1), \\
u_t(t,1)=\eta_1(t),\\
\dot\eta_1(t)=-\alpha_1\eta_1(t)-\beta_1u_x(t,1),\label{eq:eta-11}\\
u(t,0)=0,\quad t\geq 0, \label{eq:eta-0}\\
u(0, \cdot) = u_0, \quad u_t (0, \cdot) = u_1, \quad \mbox{ on } (0, 1), \\
\eta_1(0)=\eta_0.
\end{subnumcases}
where $a(\cdot)$, $q(\cdot)$ are respecting \ref{item1}-\ref{item2}, and with $\alpha_1$ and $\beta_1$ are strictly positive.

The state of this system is 
\begin{align}
  \mathcal{X}_2(t)=&[ u(t,\cdot),\, u_t(t,\cdot),\, \eta_1(t)] \in \text{Dom}(\mathcal{A}_2),
\end{align} 
where $\mathcal{A}_2$ is the unbounded operator associated with \eqref{sys_prop_2}. 
The domain is defined as 
\begin{align}
  \text{Dom} (\mathcal{A}_2)=\{ z\in X_{2,s},\, z_1(0)=0,\, z_2(1)=z_3\},
\end{align}
where $X_{2,s}$ is the space of strong solutions, and $X_{2,w}$ is the space of weak solutions defined as
\begin{eqnarray}
  X_{2,s}=H^2\times H^1 \times \mathbb{R},\\
  X_{2,w}=H^1\times L^2 \times \mathbb{R}.
\end{eqnarray}
Finally, consider 
\begin{align}
  \Gamma_2 (\mathcal{X}_2(t))=&\int_{0}^{1} (u_t^2(t,x)+u_x^2(t,x))dx \notag\\&+ \eta_1^2(t).
\end{align}
Then, there exist a positive constant $\rho$ and a positive constant $M$ 
such that for every weak solution $\mathcal{X}_2$, it holds 
\begin{align}
  \Gamma_2 (\mathcal{X}_2(t)) \leqslant M \Gamma_2 (\mathcal{X}_2(0)) e^{-\rho t}.
\end{align}
And the system is exponentially stable towards the origin of $X_{2,w}$.

In addition, it holds that $\max_{x\in[0,1]}\vert u(t,x)\vert$ tends exponentially to zero as $t$ tends to infinity, with a decay rate larger than or equal to $\rho$.
\end{proposition}

{ The following system is when we consider an integral part at the dynamic boundary for \eqref{sys_prop_2}. }

\begin{proposition}\label{prop:W2-D}
Consider the following 1D wave equation, 
\begin{subnumcases}{\label{sys_prop_3}}
u_{tt} - (au_x)_x  = -qu_t, \quad (t,x)\in \mathbb{R}^+ \times (0, 1), \\
u_t(t,1)=\eta_1(t),\\
\dot\eta_1(t)=-\alpha_1\eta_1(t)-\alpha_2\eta_2(t)-\beta_1u_x(t,1),\label{eq:xi111}\\
\dot\eta_2(t)=\eta_1(t),\label{eq:eta2}\\
u(t,0)=0,\quad t\geq 0,\label{eq:eta-00} \\
u(0, \cdot) = u_0, \quad  u_t (0, \cdot) = u_1, \quad \mbox{ on } (0, 1), \\
\eta_1(0)=\eta_0, \quad \eta_2(0)=\eta_{2,0}.
\end{subnumcases}
where $a(\cdot)$, $q(\cdot)$ are respecting \ref{item1}-\ref{item2}, and with $\alpha_1$, $\alpha_2$ and $\beta_1$ are strictly positive.

For $x\in [0,1]$, define 
\begin{align}
v(x)&:=C_2\int_0^x\frac{ds}{a(s)}, \\
C_2&:=\frac{a(1)\alpha_2}{a(1)\alpha_2\int_0^1\frac{ds}{a(s)}+\beta_1}u_*(1). 
\end{align}
The state of this system is 
\begin{align}
  \mathcal{X}_3(t)=[& u(t,\cdot),\, u_t(t,\cdot),\, \eta_1(t),\notag\\& \eta_2(t)] \in \text{Dom}(\mathcal{A}_3),
\end{align} 
where $\mathcal{A}_3$ is the unbounded operator associated with \eqref{sys_prop_3}. 
The domain is defined as 
\begin{align}
  \text{Dom} (\mathcal{A}_3)=\{ z\in X_{3,s},\, z_1(0)=0,\, z_2(1)=z_3\},
\end{align}
where $X_{3,s}$ is the space of strong solutions, and $X_{3,w}$ is the space of weak solutions defined as
\begin{align}
  X_{3,s}=H^2\times H^1 \times \mathbb{R}^2,\label{def_X3s}\\
  X_{3,w}=H^1\times L^2 \times \mathbb{R}^2.\label{def_X3w}
\end{align}
Finally, consider 
\begin{align}
  \Gamma_3 (\mathcal{X}_2(t))=&\int_{0}^{1} ((u(t,x)-v(x))^2 +u_t^2(t,x)\notag\\&+u_x^2(t,x))dx \notag\\& + \eta_1^2(t)+(\eta_2(t)-\frac{\beta_1C_2}{a(1)\alpha_2})^2.
\end{align}
Then, there exists a positive constant $\rho$ and a positive constant $M$ such that
for every weak solution $\mathcal{X}_2$, it holds 
\begin{align}
  \Gamma_3 (\mathcal{X}_3(t)) \leqslant M\Gamma_3 (\mathcal{X}_3(0)) e^{-\rho t}.
\end{align}
And the system is exponentially stable towards the attractor defined as $\text{ker}\ (\Gamma_3(\cdot))$.

In additions, it holds that $\max_{x\in[0,1]}\vert u(t,x)-v(x)\vert$ tends exponentially to zero as $t$ tends to infinity, with a decay rate larger than or equal to $\rho$.
\end{proposition}
{ Now we consider the case where $\alpha_2=0$ in \eqref{sys:propW2-W1}.
 This system has been studied in a more general and multidimensional setup in \cite{BUFFE2017207}, the author establishes with lessen hypothesis logarithmic decay rates. }
\begin{proposition}\label{prop:W1-W1}
Consider the following 1D wave equation, 
\begin{subnumcases}{\label{sys:propW1-W1}}
u_{tt} - (au_x)_x  = -qu_t, \ (t,x)\in \mathbb{R}^+ \times (0, 1), \\
u_t(t,1)=\eta_1(t),\\
u_t(t,0)=\xi_1(t),\\
\dot\eta_1(t)=-\alpha_1\eta_1(t)-\beta_1u_x(t,1),\label{eq:xi11-1}\\
\dot\xi_1(t)=-\gamma_1\xi_1(t)+\mu_1u_x(t,0),\label{eq:xi10-0}\\
u(0, \cdot) = u_0, \quad  u_t (0, \cdot) = u_1 \quad \mbox{ on } (0, 1), \\
\eta_1(0)=\eta_0, \quad  \xi_1(0)=\xi_0.
\end{subnumcases}
where $a(\cdot)$, $q(\cdot)$ are respecting \ref{item1}-\ref{item2}, and with $\alpha_1,\beta_1$,  $\gamma_1$ and $\mu_1$ are positive.  
The state of this system is 
\begin{align}
  \mathcal{X}_4(t)=&[ u(t,\cdot),\, u_t(t,\cdot),\, \eta_1(t),\, \xi_1(t)] \in \text{Dom}(\mathcal{A}_4),
\end{align} 
where $\mathcal{A}_4$ is the unbounded operator associated with \eqref{sys:propW1-W1}. 
The domain is defined as 
\begin{align}
  \text{Dom} (\mathcal{A}_4)=\{& z\in X_{3,s};\notag\\& z_2(0)=z_4,\, z_2(1)=z_3\},
\end{align}
where $X_{3,s}$ is the space of strong solutions, and $X_{3,w}$ is the space of weak solutions, both defined in \eqref{def_X3s}-\eqref{def_X3w}
Finally, consider 
\begin{align}
  \Gamma_4 (\mathcal{X}_4(t))=&\int_{0}^{1} (u_t^2(t,x)+u_x^2(t,x))dx \notag\\&+ \eta_1^2(t) + \xi_1^2(t).
\end{align}
Then, there exists a positive constant $\rho$ and a positive constant $M$ 
such that, for every weak solution $\mathcal{X}_4$, it holds 
\begin{align}
  \Gamma_4 (\mathcal{X}_4(t)) \leqslant M \Gamma_4 (\mathcal{X}_4(0)) e^{-\rho t}.
\end{align}
And the system is exponentially stable towards the attractor $S_4$ defined by
\begin{align}
  S_4 =\{&z\in X_{3,w},\ z_1(\cdot)=d,d\in\mathbb{R},\, z_2(\cdot)=0,\notag\\ & z_3 = 0,\, z_4=0\}.
\end{align}
which is the kernel of $\Gamma_4(\cdot)$.

In addition, there exists $u_*$ so that $\max_{x\in[0,1]}\vert u(t,x)-u_*\vert$ tends exponentially to zero as $t$ tends to infinity.
\end{proposition}
\begin{pf}

We start by proving Proposition~\ref{prop:W1-D}. As before, the argument is based on an appropriate Lyapunov function $\bar{V}_u=\bar{F}_u+\ell \bar{W}_u$
where $\ell$ is a positive constant  to be chosen and
{
\begin{align}
\bar{F}_u(t):=&\frac12\int_0^1(u_t^2+au_x^2)\, dx+\frac{a(1)}{2\beta_1}\eta_1^2,  \label{Fbar}\\
\bar{W}_u(t):=&\int_0^1uu_t\, dx+\frac12\int_0^1qu^2\, dx\notag\\&+\frac{a(1)}{\beta_1}\eta_1u(t,1)
{
+\frac{a(1)\alpha_1}{2\beta_1}u(t,1)^2
}
.  
\end{align}
One gets, using integration by parts, \eqref{eq:eta-11}, and \eqref{eq:eta-0}
\begin{align}
\dot{\bar{F}}_u(t):=&-\int_0^1(qu_t^2)\, dx-\frac{a(1)\alpha_1}{\beta_1}\eta_1^2,  \\
\dot{\bar{W}}_u(t):=&\int_0^1u_t^2\, dx-\int_0^1au_x^2\, dx+\frac{a(1)}{\beta_1}\eta_1^2.  
\end{align}
}
Therefore
\begin{align}
\dot {\bar{V}}_u=&-\int_0^1(q-\ell)u_t^2\, dx-\ell\int_0^1au_x^2\, dx\notag\\& -\frac{a(1)}{\beta_1}(\alpha_1-\ell)\eta_1^2.
\end{align}
The conclusion follows by taking $\ell>0$ small enough and noting that, thanks to the Dirichlet boundary condition \eqref{eq:eta-00}, for every $t\geq 0$ and $x\in [0,1]$
\begin{align}
\vert u(t,x)\vert^2&=\vert u(t,x)-u(t,0)\vert^2\notag\\& \leq \int_0^1u_x^2(t,x)\, dx\leq 2E_u(t)\leq 2\bar{F}_u(t).   \label{eq:final-1}
\end{align}
One proceeds by establishing an analog to Proposition~\ref{prop1} where $\Gamma$ and $V$ are replaced by $\Gamma_2$ and $\bar{V}_u$ in order first to obtain that $\dot{\bar{V}}_u\leq -\rho \bar{V}_u$ for some positive constant $\rho$ independent of the state and finally to conclude as in the final part of the argument of Theorem~\ref{th:main1}.
  
We next turn to the proof of Proposition~\ref{prop:W2-D}. Using the notations of the proposition, we set 
\begin{align}
w(t,x):&=u(t,x)-v(x),\quad t\geq 0,\ x\in [0,1], \notag\\  \bar{\eta}_2(t):&=\eta_2(t)+\frac{\beta_1C_2}{a(1)\alpha_2},\quad t\geq 0 \label{def_eta_bar}.  
\end{align}
It is a matter of elementary computations to check that {$w$ is the solution of \eqref{sys_prop_3} with different and corresponding initial conditions}
with the Dirichlet boundary condition at $x=0$ (since $v(0)=0$) and the boundary condition given by 
\begin{align}
\dot\eta_1(t)=&-\alpha_1\eta_1(t)-\alpha_2\bar{\eta}_2(t)-\beta_1w_x(t,1),\label{eq:xi1-w}\\
\dot{\bar{\eta}}_2(t)=&\eta_1(t).\label{eq:bareta2}
\end{align}
{
It holds
\begin{align}
w_*(1)&:= w(t,1)-\bar{\eta}_2(t)\notag\\&=u_*(1)-v(1)-\frac{\beta_1C_2}{a(1)\alpha_2}\notag\\&=u_*(1)-C_2(\int_0^1\frac{ds}{a(s)}+\frac{\beta_1C_2}{a(1)\alpha_2})=0.
\end{align}}
We have essentially reduced the problem to only deal with solutions of \eqref{sys:propW1-W1} with the Dirichlet boundary condition at $x=0$, with the additional constraint that $w_*(1)=0$.  In that case, we consider the candidate Lyapunov function 
$\tilde{V}_w=\tilde{F}_w+\ell \tilde{W}_w$ where $\ell$ is a positive constant to be chosen and 
\begin{align}
{\tilde{F}_w(t):}&={\frac12\int_0^1(w_t^2+aw_x^2)\, dx}\notag\\ &+\frac{a(1)}{2\beta_1}(\eta_1^2+\alpha_2\bar{\eta}_2^2),\label{F_tilde_bar}\\
\tilde{W}_w(t):=&\int_0^1ww_t\, dx+\frac12\int_0^1qw^2\,  dx\notag\\&+\frac{a(1)}{\beta_1}\Big(\frac{\alpha_1}2\bar{\eta}_2^2+\eta_1\bar{\eta}_2\Big).  
\end{align}
One gets
\begin{align}
\dot {\tilde{V}}_u=&-\int_0^1(q-\ell)w_t^2\, dx-\ell\int_0^1aw_x^2\, dx\notag\\&-\frac{a(1)}{\beta_1}(\alpha_1-\ell)\eta_1^2-\ell \frac{a(1)\alpha_2}{\beta_1}\bar{\eta}_2^2,  
\end{align} 
where we have repeatedly used the equality $w(t,1)= \bar{\eta}_2(t)$. By following what has been done previously, the conclusion follows. 
  
We finally prove Proposition~\ref{prop:W1-W1}. As before the argument is based on an appropriate Lyapunov function $\bar{V}_u$ defined later.
We first consider $F$ given in \eqref{eq:E1} and note that for $t\geq 0$ it holds
\begin{equation}\label{eq:w1w1-Fu}
\dot F=-\int_0^1qu_t^2\, dx-\frac{a(1)\alpha_1}{\beta_1}\eta_1^2(t)-\frac{a(0)\gamma_1}{\mu_1}\xi_1^2(t).
\end{equation} 
We next compute along solutions of \eqref{sys:propW1-W1} the following time derivative
{
\begin{align}  
\frac{d}{dt}&\Big(\int_0^1\Big(u(t,x)-u(t,1)\Big)u_t(t,x)\,dx\Big)=\notag\\ &+\int_0^1u_t^2\, dx-\int_0^1q\big(u(t,x)-u(t,1)\big)u_t(t,x)\, dx\notag\\
&+\Big(u(t,1)-u(t,0)\Big)a(0)u_x(t,0)\notag\\&-\int_0^1au_x^2\, dx
{-\eta_1\int_0^1u_t(t,x)\, dx}. \label{eq_uu1}
\end{align}}
In the above equation, we use \eqref{eq:xi10-0} to get rid of $u_x(t,0)$ and, to obtain for $t\geq 0$ that
\begin{align}  
&\Big(u(t,1)-u(t,0)\Big)u_x(t,0)=\notag\\ & +\Big(u(t,1)-u(t,0)\Big)\frac{\dot\xi_1+\gamma_1\xi_1}{\mu_1}
\notag\\&=\frac{d}{dt}\Big(\Big(u(t,1)-u(t,0)\Big)\frac{\xi_1}{\mu_1}\Big)\notag\\
\ &-(\eta_1-\xi_1)\frac{\xi_1}{\mu_1}+\frac{\gamma_1\xi_1}{\mu_1}\Big(u(t,1)-u(t,0)\Big).
\end{align}

Setting for $t\geq 0$
\begin{align}
G_u(t):&=\int_0^1\Big(u(t,x)-u(t,1)\Big)u_t(t,x)\,dx\notag\\&-\frac{a(0)\xi_1}{\mu_1}\Big(u(t,1)-u(t,0)\Big),  
\end{align}
we deduce from the above that along with solutions of \eqref{sys:propW1-W1} that
{
\begin{align}  
\dot G_u&=\int_0^1u_t^2\, dx-\int_0^1au_x^2\, dx
-\frac{a(0)\xi_1}{\mu_1}(\eta_1-\xi_1)\notag\\&+\frac{a(0)\gamma_1\xi_1}{\mu_1}\Big(u(t,1)-u(t,0)\Big)
{-\eta_1\int_0^1u_t(t,x)\, dx}
\notag\\& -\int_0^1q\Big(u(t,x)-u(t,1)\Big)u_t\, dx
.\label{eq:gu}
\end{align}}
We finally recall that there exists a positive constant $C_a$ (independent of the solutions of \eqref{sys:propW1-W1}) such that, for $t\geq 0$,
\begin{align}\label{eq:astuce}
\int_0^1\vert u_t\vert\,dx+&\max_{x\in [0,1]}\vert u(t,x)-u(t,1)\vert\notag\\ &\leq  \int_0^1(\vert u_t\vert+
\vert u_x\vert)\,dx\notag\\&\leq C_aE_u^{1/2}(t).
\end{align}
We now choose $\bar{V}_u=F+\ell G_u$ for $\ell>0$ small enough. Using repeatedly the Cauchy-Schwarz inequality, and \eqref{eq:astuce} in \eqref{eq:w1w1-Fu} and \eqref{eq:gu}, one gets for $\varepsilon$ and $\ell$ small enough that 
\eqref{eq:ET-est0} and \eqref{eq:lyap} hold true, from which one deduces Item$(i)$ of Proposition~\ref{prop:W1-W1}.
  
Finally, to get Item$(ii)$ of Proposition~\ref{prop:W1-W1}, first notice that $u(t,1)$ admits a limit $u_*$ as $t$ tends to infinity since, for every $t,t'>0$ it holds $u(t,1)-u(t',1)=\int_{t'}^t\eta_1$ and $\eta_1$ decreases to zero exponentially. The conclusion follows now by using \eqref{eq:astuce}.
\end{pf}

\begin{remark}\label{rem:common-trick}
In the proofs of all our results, one could use the function $G_u$ (especially the integral term) to obtain the exponential decrease of $E_u$ and some of the components of the Wentzell's boundary conditions. However, this does not allow one to determine the limit $u_*$ for the solution $u$ in terms of initial conditions. In particular, we are not able to characterize
$u_*$ in Proposition~\ref{prop:W1-W1}. 

Note also that 
\begin{align}
  u(t,x)-u(t,1)=-\int_x^1 u_x(t,s)ds.
\end{align}
This can be related with the means of $u_x$ and therefore we have extended our Lyapunov function with a space moving evaluation of the mean of the force/torque.
Indeed $u_x$ is associated with the torque or the force in mechanical setup. 

\end{remark}

\section{Numerical schemes and simulations.}

There exist several ways to compute numerical approximation of the solution of evolution problems associated with partial differential equation, \cite{tadmor2012review}. In the case under consideration, spectral methods lead to an estimation of the base function at each time step due to the dynamics boundary condition. This requires an important computing power. As we have only one dimension in space finite-element methods reduce to finite difference methods with (possibly unequal) spacial step. Finite different methods can be delicate to design in order to ensure at the same time numerical stability and good approximation. Note that there also exist specific schemes based on Riemann invariants \cite{boldo2013wave}. These last schemes have good numerical property, but their extension to dynamic boundary conditions is not obvious.

In this paper, we suggest a new approach, which provides numerical scheme stability and therefore achieves structural stability. It is based on the discretization of the Lagrangian associated with the wave equation. This approach leads to a special finite difference scheme. As previously said the wave equation in its stationary form can be associated with a Lagrangian. For the case under consideration \eqref{sys_problem}, (in the stationary case where $\nu=U(t)=f_c=\gamma_1=f_{ac}=0$), this Lagrangian is given by
\begin{align}
L(v(t,\cdot))&=\int_0^1 \frac{1}{2}(v_t^2(t,x)-a(x)v_x^2(t,x))dx\notag\\&+\frac{1}{2}(\frac{a(1)}{\beta_1}{v_t(t,1)}^2+\frac{a(0)}{\mu_1}{v_t(t,0)}^2) \label{eq_lagragien}.
\end{align}
Following the strategy in \cite{kot2014first} and the least action principle, the dynamics of the system is associated with a stationary action. The action for any time interval is given as
\begin{align}
I(v)=\int_{t_i}^{t_f} L(v(t,\cdot))dt.
\end{align}
A stationary action means that the first variation is equal to zero
\begin{align}
\delta I(v,\delta v) = 0,
\end{align}
where the first variation is defined as
\begin{align}
\delta I (v,\delta v)= \delta I(v+\delta v) - \delta I(v) + O(\norme{\delta v}^2).
\end{align}
Computation gives the following stationary system
\begin{subnumcases}{}
v_{tt} - (av_x)_x  =0,\quad (x,t) \mbox{ in } \mathbb{R}^+ \times (0, 1), \\
v_{tt}(t,1)=-\beta_1v_x(t,1),\\
 v_{tt}(t,0)=\mu_1v_x(t,0).
\end{subnumcases}
{ This is the stationary part of \eqref{sys_problem}, as usual the less action principle, the dissipation and the input are added afterward to obtain exactly \eqref{sys_problem}.}
Now consider a discrete version of \eqref{eq_lagragien}
\begin{align}
L_d(v_d(t)[\cdot])&=\frac{1}{2}\sum_{i=1}^{N-1} [\dot{v}_d(t)[i]^2  \notag\\&-\frac{a_{i-1}}{12}(\frac{v_d(t)[i]-v_d(t)[i-1]}{dx_i})^2\notag \\ 
-& \frac{a_{i}}{3}(\frac{v_d(t)[i+1]-v_d(t)[i-1]}{dx_i+dx_{i+1}})^2\notag\\&-\frac{a_{i+1}}{12}(\frac{v_d(t)[i+1]-v_d(t)[i]}{dx_{i+1}})^2]dx_{i} \notag\\
+& \frac{1}{2}\frac{a_N}{\beta_1} \dot{v}_d(t)[N]^2+\frac{1}{2}\frac{a_0}{\mu_1} \dot{v}_d(t)[0]^2.
\end{align}
The integral part in $v_x^2$ has approximated using Simpson's $1/3$ rule. The derivation of the Euler-Lagrange equation can then be done by a symbolic numerical computation. This gives an autonomous stationary linear finite dimensional system: 
\begin{align}
E\ddot{v}_d(t)=Av_d(t),
\end{align} 
with $\sigma(A)\in {\bf i} \mathbb{R}$. It holds
\begin{align}
E=\text{diag}\begin{bmatrix}
  \frac{1}{\beta_1} & dx_1 & \hdots & dx_{N-1} & \frac{1}{\mu_1} 
\end{bmatrix}.
\end{align}
Then we add dissipation with a positive symmetric matrix $R$, source term (disturbance and action) and observation,
\begin{subnumcases}{}
E\ddot{v}_d(t)=Av_d(t) -R\dot v_d(t)+BU(t)+f_{et}, \\
y(t)=C \dot v_d(t).
\end{subnumcases}
{with 
\begin{align}
  f_{et}^T=\begin{bmatrix}
    f_{c}&f_1&f_2&\hdots & f_{ac}
\end{bmatrix}
\end{align}
 which represents the disturbance, and with
 \begin{align}
  R=\text{diag}\begin{bmatrix}
    \nu & q_1 & \hdots & q_{N-1} & \gamma_1 
  \end{bmatrix}.
  \end{align}}
The control $U(t)$ is computed through
\begin{subnumcases}{}
{ \dot\eta_v(t)} = y(t)-y_\text{ref},\\
U(t)=-k_i { \eta_v(t)}  -k_p (y(t)-y_\text{ref}).
\end{subnumcases}
As the main idea of this discretization scheme is to have a good approximation of the energy, we suggest going on with this idea using symplectic integrator scheme, see \cite{Donnelly2005} and references within. These methods, like the Crank-Nicolson method have the property preserve the energy as time evolves. It is known that for a system which has an eigenvalue in ${\bf i} \mathbb{R}$ explicit schemes are unstable, and implicit schemes are exponentially stable see \cite{Donnelly2005}. As our system has structurally the zero eigenvalue, symplectic numerical discretization schemes tend to give better behaviors approximation.

The idea of a symplectic scheme is to combine an implicit scheme together with an explicit one. This leads to 
\begin{align}
v_d[k+1]=&v_d[k]+\Delta t\,\dot{v}_d[k+1],\\
E\dot v_d[k+1] =& E\dot v_d[k] + \Delta t\, Av_d[k]-\Delta t\, R\dot v_d[k+1]\notag\\& +\Delta t\, BU[k]+\Delta t f.
\end{align} 
The second line is implicit, but $R$ in our case is a diagonal matrix and so the associated inverse matrix is easily computed
\begin{align}
\dot v_d[k+1] =& (1+ \Delta t E^{-1} R)^{-1} (\dot v_d[k]+\Delta t\, E^{-1}Av_d[k]\notag\\&+\Delta t\, E^{-1}BU[k]+\Delta t E^{-1} f).
\end{align}
There are several key points to note in this last equation. First, the term $(1+E^{-1}\Delta t R)^{-1}$ correspond to a contraction map in the case where $R$ is positive, and therefore is associated with dissipation terms. Second, in the case where $R$ represent anti-dissipation term, there exist discretized steps $\frac{\Delta t}{dx}$ where the numerical shame is undefined. Third, where $R=0$, these equations are two-step explicit ones. 
The value selected for the numerical simulation for the output regulation problem is summarized in Table~\ref{Table_simu}.

\begin{table}[h!]
\centering
\begin{tabular}{c   c | c c}
Symbol &   value & Symbol & value  \\
\hline 
\hline
$N$&    $199$& $f_c$ & $-1$  \\[1.05ex]
$a(x)$  & $\sin(2x)+2 $& $f_{ac}$ & $1$  \\[1.05ex]
$q(x)$  & $.01+.1x^2$& $k_p$  & $10$ \\[1.05ex]
$f(x)$  & $sin(2\pi x)$& $k_i=\alpha_2$  & $20$ \\[1.05ex]
$\beta_1$ & $20$ & $v_1^{\text{ref}}$  & $5$ \\[1.05ex]
$\mu_1$ & $20$&  $v_d[\cdot]$  & $0$ \\[1.05ex]
$\nu $ & $1$ &  $\dot v_d[0:N]$  & $0$ \\[1.05ex]
$\gamma_1$ & $1$ &  $\Delta t$  & $0.001$ 
\end{tabular}
\caption{Parameter values for the simulation. \label{Table_simu}}
\end{table}

\begin{figure}[!ht]
  \centering
  \includegraphics[width=.9\linewidth]{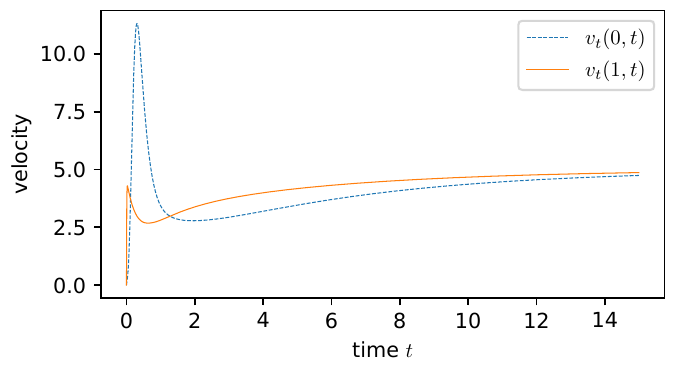}
  \caption{The boundary velocities times responses.\label{fig_B}}
\end{figure}

\begin{figure}[!ht]
  \centering
  \includegraphics[width=.9\linewidth]{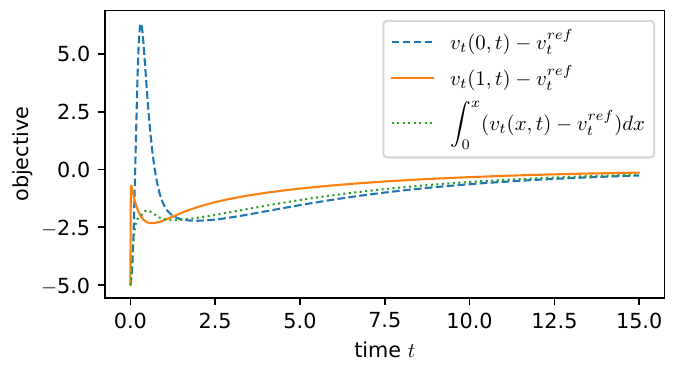}
  \caption{The objectives times responses.\label{fig_obj}}
\end{figure}

\begin{figure}[!ht]
  \centering
  \includegraphics[width=.9\linewidth]{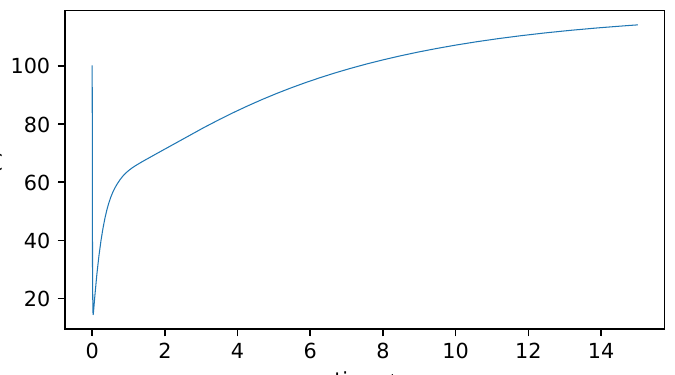}
  \caption{The control law time response.\label{fig_U}}
\end{figure}

\begin{figure}[!ht]
  \centering
  \includegraphics[width=.9\linewidth]{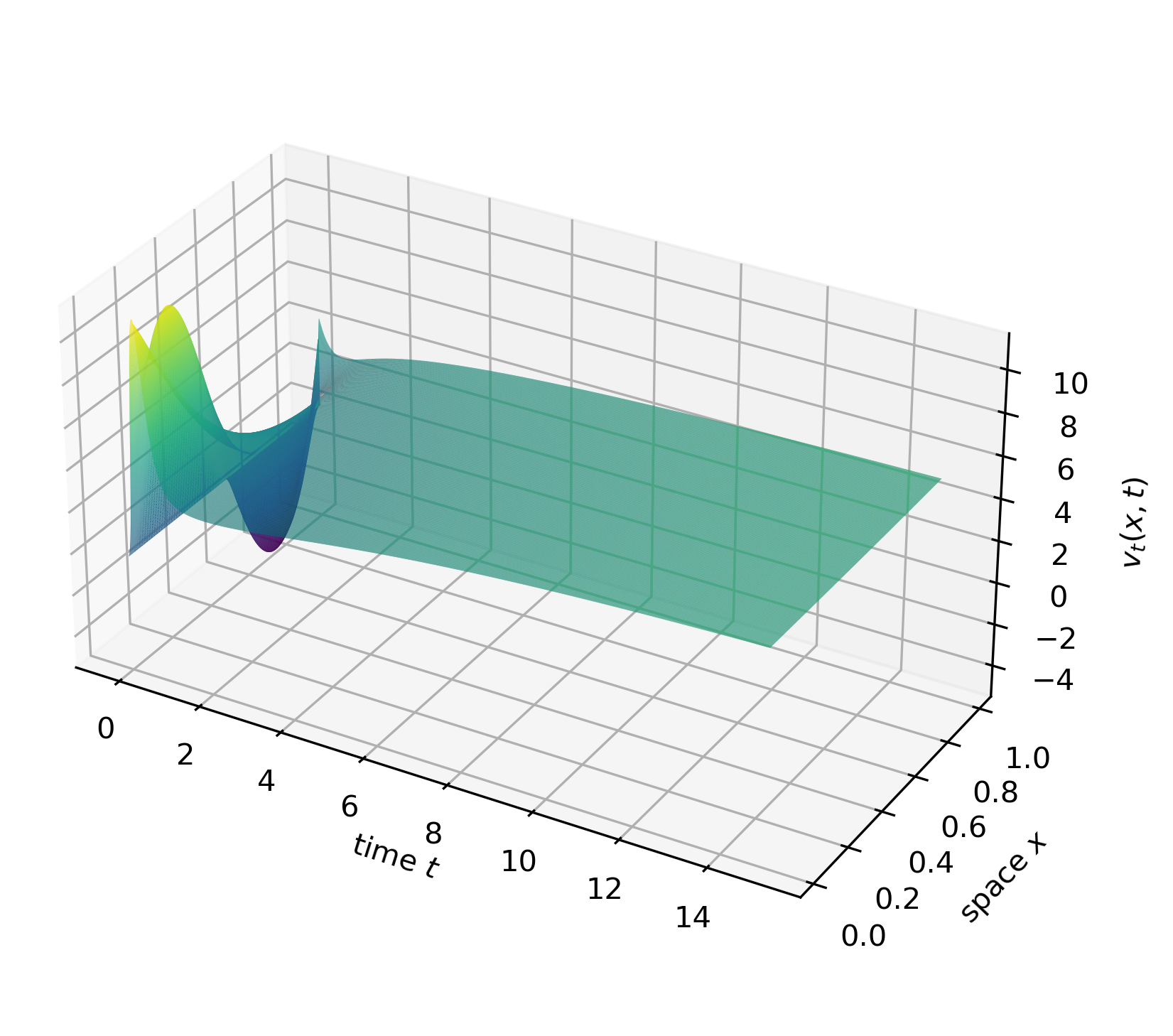}
  \caption{The distributed velocity $\dot{v}(t,x)$ time response.\label{fig_ut}}
\end{figure}

\begin{figure}[!ht]
  \centering
  \includegraphics[width=.9\linewidth]{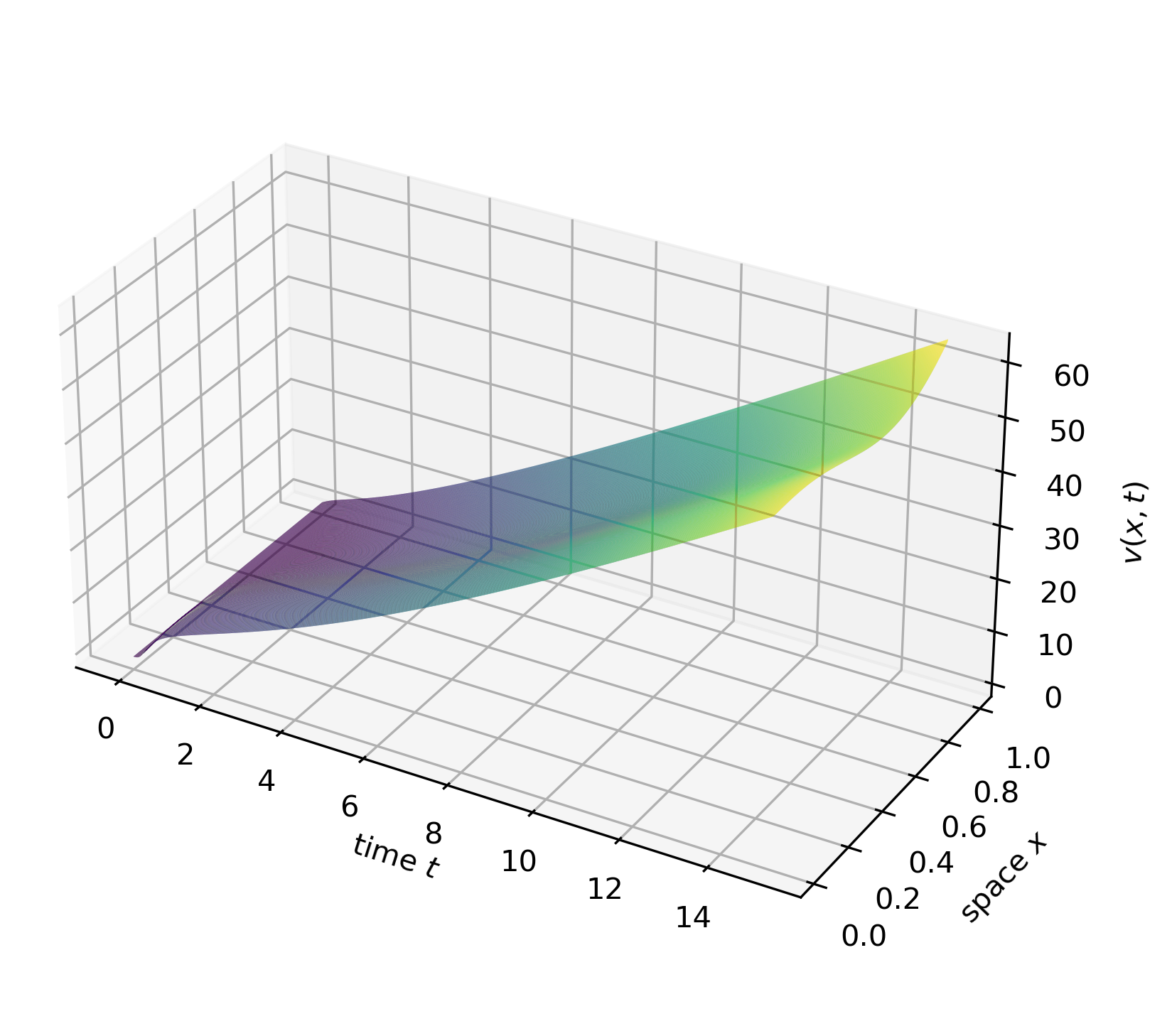}
  \caption{The distributed position $v(t,x)$ time response.\label{fig_u}}
\end{figure}

\begin{figure}[!ht]
  \centering
  \includegraphics[width=.9\linewidth]{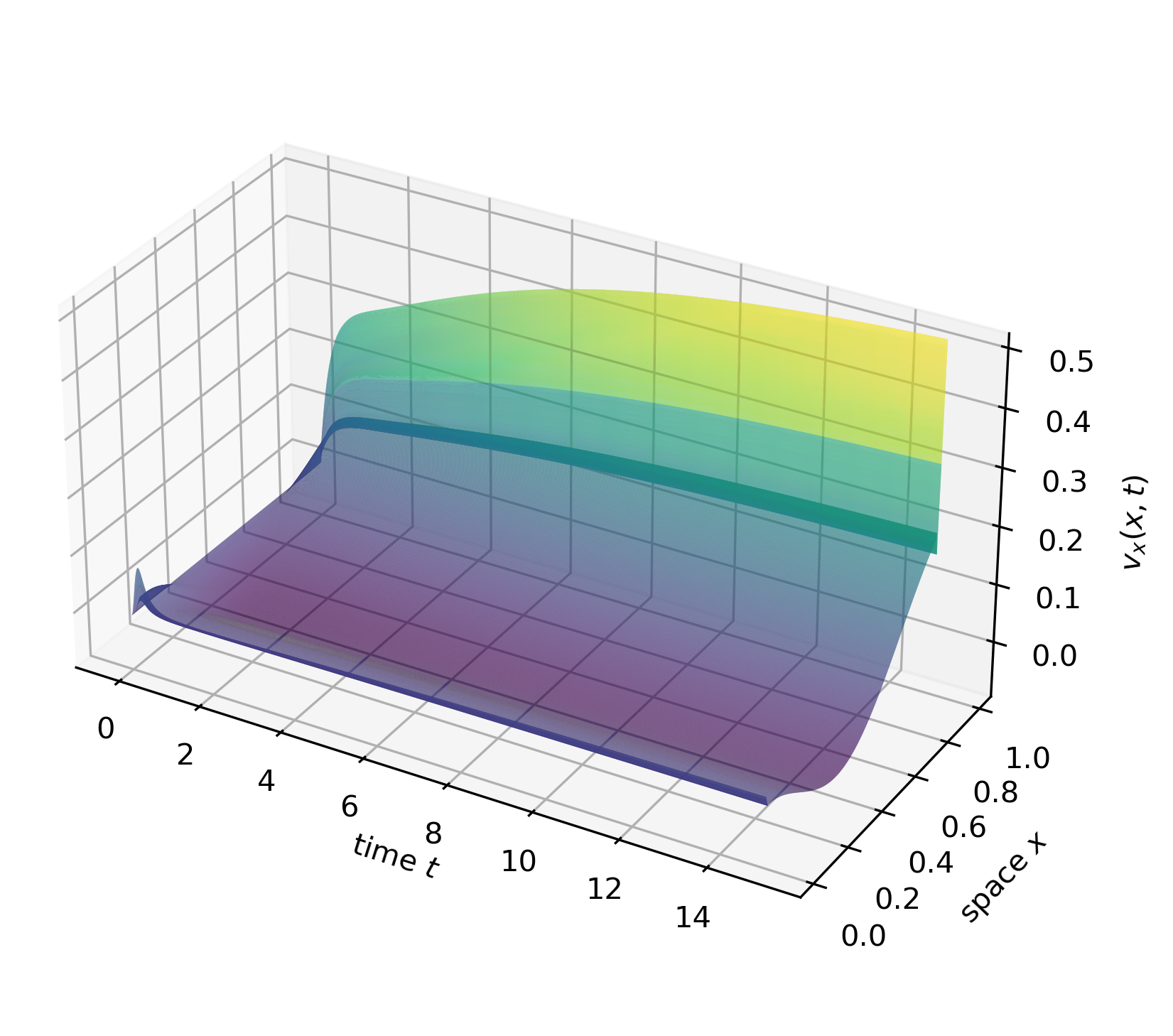}
  \caption{The distributed $v_x(t,x)$ time response.\label{fig_ux}}
\end{figure}

{
The Figure~\ref{fig_B} illustrates the behavior of the output regulation problem, we observe that boundary velocities of the system goes exponentially towards the constant reference. In Figure~\ref{fig_obj} the time response of the regulation problem objectives are depicted. The in-domain velocity converges in $L^2$ norm towards the reference. The control law associated with these time responses are given in Figure~\ref{fig_U}. It is not clear how to select the control gain to provide rapidity and robustness. Getting an urge integrator gain in order to have the control go faster towards its steady state may cause some heavy oscillation. However, as proven the exponential stability still holds.\par

The time response of the wave equation velocity is drawn as a surface in a 3d perspective in Figure~\ref{fig_ut}. There is first some important oscillation, with traveling wave going back and forth from the boundary, then the oscillation rapidly goes smaller, and finally the velocity goes smoothly towards the reference. The time response of the position is given in Figure~\ref{fig_u}. The impact of the constant disturbance are more visible in this graph.  The oscillations observed in Figure~\ref{fig_u} are mainly due to the disturbance which needs a particular distribution of the position along the space. Once this particular distribution is obtained, the constant disturbance is compensated by the integrator. The last figure, Figure~\ref{fig_ux} depicts $v_x(t,\cdot)$ it allows to observe the effect of the disturbance and the in-domain damping. The smooth convergence of the velocity can be compared with the behavior of $v_x(t,\cdot)$.
}

\section{Conclusion}
{ This paper presents the first systematic  Lypunov analysis
for a $1$-dimensional damped wave equation subject to various dynamic (or Wentzell) boundary conditions, in the case where the damping is everywhere active. As a particular case, we also provide a regulation law for a wave equation (with dynamic boundary conditions) by the means of a PI control.
The control law achieved exponential decay rate towards the constant reference, and the rejection of constant disturbance.
The possible rejection of the disturbance by the integral action can be explained by the interne model principle. The numerical simulation shows the behavior of the closed-loop system with unknown disturbance.
Future work will be to use some of the exponential decay system study in the appendix as the target system for infinite-dimensional backstepping control design.
There is also a great interest towards considering non-linear terms. For example, what is happening when the damping is non-linear like in \cite{GUGAT201472}, or even can we generalize towards non-linear waves as
\begin{align}
  u_{tt}=(\frac{a(x)u_x}{\sqrt{1+u_x^2}})_x-q(x)u_t.
\end{align}
Moreover for practical applications there is great interest studying the wave equation with dynamics boundary condition but with a non-linear friction term at the boundary opposite to the control, typically LuGre friction term.
}

\bibliographystyle{plain}       
\bibliography{bibli}  

\appendix

\section{Proof of Theorem \ref{th:wellposeness} \label{app_well_posed}}
The proof follows the same lines as the ones exposed in \cite{roman2021}. The idea of the proof is to decompose the operator $\mathcal{A}$ defined in \eqref{def_A} into a maximal monotone part and a remaining part. We should be able to cancel the remaining part with a bijective change of variable. Finally, we conclude using the following theorem.
\begin{theorem}[{Hille-Yosida \cite[Theorem~7.4 ]{brezis2010functional}}]\label{theorem_Hille-Yosida_inhomo}

Let $\mathcal{A}$ be a maximal operator on the Hilbert space $H$ then for every $X_0\in \mathcal{D} \big( \mathcal{A}\big)$ there exists a unique solution $X$ to the following abstract problem.
\begin{subnumcases}{\label{BTU:def_abstract_problem_inhomo}}
\frac{dX}{dt}(t)+\mathcal{A}X(t)=0, \\
X(0)=X_0.
\end{subnumcases}
with
\begin{equation}
X\in C^1([0,\infty); H) \cap C([0,\infty); \mathcal{D}\big(\mathcal{A}\big) ).
\end{equation}
\end{theorem}

Now consider the following operator
\begin{equation}
\forall z\in \mathcal{D}(G), \ Gz= \begin{bmatrix}
-z_2 \\[1ex] -(az'_1)'+z_2+z_1 \\[1ex]  \beta_1\,z'_1(1) \\[1ex] 0 \\[1ex]  -\mu_1\,z'_1(0)
\end{bmatrix}\label{def:op_G},
\end{equation}
and the following matrix
\begin{equation}
B = \begin{bmatrix}
0&0&0&0&0\\ 1&-q+1 &0&0&0 \\0&0&- \alpha_1&-\alpha_2&0\\0&0&1&0&0 \\0&0&0&0&-\gamma_1
\end{bmatrix}\label{def:op_B}.
\end{equation}
The domain of $G$ is equal to the domain of $\mathcal{A}$. One gets
\begin{equation}
\mathcal{A}=G+B.
\end{equation}
$G$ is a monotone part, this is established in the following lemma and $B$ is a bounded operator. 
\begin{lemma}\label{BTU:lem_maximal_monotone} The unbounded linear operator $G$ defined in \eqref{def:op_G} is a maximal monotone operator on $X_w$ defined in \eqref{def:SS}.
\end{lemma}
\begin{pf}
Considering the following scalar product on $X_w$ 
\begin{align}
\scal{z}{q}=&\int_0^1 (z_1\nu +z_2q_2+az_1'\nu ') dx +\notag\\& \frac{a(1)}{\beta_1}z_3q_3+z_4q_4+\frac{a(0)}{\mu_1}z_5q_5,\\
\scal{z}{Gz}&= \int_0^1[-z_1z_2+z_2(-(az_1')'+z_2+z_1)\notag\\&-a(x)z'_1z'_2]dx+a(1)z_3z'_1(1)\notag\\&-a(0)z_5z'_1(0),
\end{align}
using integration by parts and the fact that $z\in\mathcal{D}(\mathcal{A})$, one obtains
\begin{align}
\scal{z}{Gz}= \int_0^1z_2^2 dx  \geqslant 0.
\end{align}
Thus the operator $G$ is monotone (see \cite[Chapter~7 on Page~181]{brezis2010functional}) on the Hilbert $X_w$. In addition, if we establish that
\begin{equation}
R(I+G)=X_w,
\end{equation}
then the operator $G$ is maximal monotone (see \cite[Chapter~7 on Page~181]{brezis2010functional}), $R$ stands for the range of the operator.
Let $y\in X_w$, we have to solve
\begin{equation}
z \in\mathcal{D}(\mathcal{A}),\quad z+Gz=y,
\end{equation}
which means that
\begin{align}
z_1-z_2=&y_1 \label{BTU:eq_wp_set1}, \\
z_2-(az'_1)'+z_2+z_1=&y_2 ,\\
z_3+\beta z_1'(1)=&y_3 ,\\
z_4+0=&y_4,\\
z_5-\mu_1z'_1(0)=&y_5, \label{BTU:eq_wp_set5}
\end{align}
using the fact that $z\in\mathcal{D}(\mathcal{A})$ one gets
\begin{align}
3z_1-(az_1')'&=2y_1+y_2 \label{BTU:eq_wp_st1},\\
\beta_1z'_1(1)+z_1(1)&=(y_3+y_1(1)),\\
-\mu_1z'_1(0)+z_1(0)&=(y_5+y_1(0)) \label{BTU:eq_wp_st3}.
\end{align}
This is a classical stationary problem (e.g., see \cite{brezis2010functional}) with Robin's boundaries conditions, using standard result (as done in \cite[Example 6, On Page 226]{brezis2010functional} ) one gets that as $2y_1+y_2 \in L^2(0,1)$,  \eqref{BTU:eq_wp_st1}-\eqref{BTU:eq_wp_st3} has a unique solution $z_1\in H^2(0,1)$. Now one can check that the element $z=(z_1,\ z_2,\ z_3,\ z_4,\ z_5)$ with
\begin{subnumcases}{\label{BTU:eq_wp_re}}
z_1 \text{ is the solution to \eqref{BTU:eq_wp_st1}-\eqref{BTU:eq_wp_st3}},\\
z_2=z_1-y_1\label{BTU:eq_wp_re1},\\
z_3=y_3-a(1)z'_1(1),\\
z_4=y_4
z_4=y_5+a(0)z_1'(0),
\end{subnumcases}
satisfies \eqref{BTU:eq_wp_set1}-\eqref{BTU:eq_wp_set5}. Moreover using \eqref{BTU:eq_wp_st1}-\eqref{BTU:eq_wp_st3} on \eqref{BTU:eq_wp_re} one gets that $z$ satisfying \eqref{BTU:eq_wp_re} is in $\mathcal{D}(\mathcal{A})$.
\end{pf}

Now, we are ready to state the proof of the well posedness of \eqref{sys:abs}. Note that the fact that $G$ is maximal monotone implies that $\mathcal{D}(\mathcal{A})$ is dense in $X_w$ (i.e., $\overline{\mathcal{D}(\mathcal{A})}=X_w$).

Using the bijective change of variable
\begin{equation}
z_e(t)=z(t)e^{Bt}\label{BTU:def_change_of_varaible},
\end{equation}
$z$ is the solution to \eqref{sys:abs} is equivalent to, $z_e\in \mathcal{D}\big(\mathcal{A}\big)$ is the solution to
\begin{subnumcases}{\label{BTU:def_abstract_problem_inhomo_G}}
\frac{d}{dt}z_e(t)+Gz_e(t)=0, \\
z_e(0)=z_0,
\end{subnumcases}
where $B$ is defined in \eqref{def:op_B} and $G$ is defined in \eqref{def:op_G}.

From Lemma \ref{BTU:lem_maximal_monotone}, using Theorem \ref{theorem_Hille-Yosida_inhomo} on \eqref{BTU:def_abstract_problem_inhomo_G}, and the change of variable \eqref{BTU:def_change_of_varaible}, one establishes (i). Using argument of density of $\mathcal{D}(\mathcal{A})$ in $X_w$, and $C_0$-semigroup theory one obtains the regularity of weak solutions.

Note that we refer the reader to \cite{kato2013perturbation}, \cite{pazy2012semigroups} for the notion weak solutions. Moreover part of the proof are inspired from \cite{conrad1998strong} and \cite{dandrea1992control} which in turn originates from \cite{Slemrod1989}.

\begin{extended}
\section{Additional materials}
    This section pertains to additional materials that are not included in the accepted version of the paper and includes links to online resources.

    The first line of $A$ is
\begin{align}
  \left[\begin{matrix}- \frac{a_{0}}{6 dx_{1}} - \frac{2 a_{1} dx_{1}}{3 \left(dx_{1} + dx_{2}\right)^{2}}, & \frac{a_{0}}{6 dx_{1}} ,& \frac{2 a_{1} dx_{1}}{3 \left(dx_{1} + dx_{2}\right)^{2}} ,& \hdots \end{matrix}\right]
\end{align}
The second line is
\begin{align}
  \bigg[\begin{matrix}\frac{a_{0}}{6 dx_{1}}, & - \frac{a_{0}}{6 dx_{1}} - \frac{a_{1}}{6 dx_{2}} - \frac{a_{2} dx_{1}}{6 dx_{2}^{2}} - \frac{2 a_{2} dx_{2}}{3 \left(dx_{2} + dx_{3}\right)^{2}},\end{matrix} \notag \\  
   \begin{matrix} & \frac{a_{1}}{6 dx_{2}} + \frac{a_{2} dx_{1}}{6 dx_{2}^{2}} ,& \frac{2 a_{2} dx_{2}}{3 \left(dx_{2} + dx_{3}\right)^{2}} ,& \hdots\end{matrix}\bigg]
\end{align}
The $i$-line for $i\in[3,N-2]$ for column $i-2$ at $i+2$ is 
\begin{align}
  \bigg[\begin{matrix}\frac{2 a_{i-2} dx_{i-2}}{3 \left(dx_{i-2} + dx_{i-1}\right)^{2}}, &   \frac{a_{i-1} dx_{i-2}+a_{i-2} dx_{i-1}}{6 dx_{i-1}^{2}}, \end{matrix} \notag \\ 
  \begin{matrix} & - \frac{2 a_{i-2} dx_{i-2}}{3 \left(dx_{i-2} + dx_{i-1}\right)^{2}} - \frac{a_{i-1} dx_{i-2}+a_{i-2}dx_{i-1}}{6 dx_{i-1}^{2}} \end{matrix} \notag \\  
   \begin{matrix}  - \frac{a_{i} dx_{i-1}+a_{i-1}dx_{i}}{6 dx_{i}^{2}} - \frac{2 a_{i} dx_{i}}{3 \left(dx_{i} + dx_{i+1}\right)^{2}}, \end{matrix} \notag \\   \begin{matrix} 
     & + \frac{a_{i} dx_{i-1}+a_{i-1}dx_{i}}{6 dx_{i}^{2}}, & \frac{2 a_{i} dx_{i}}{3 \left(dx_{i} + dx_{i+1}\right)^{2}}\end{matrix}\bigg]
  \end{align}
and zero elsewhere.
The $N-1$ line
\begin{align}
\bigg[\begin{matrix}\hdots, & \frac{2 a_{N-2} dx_{N-2}}{3 \left(dx_{N-2} + dx_{N-1}\right)^{2}} ,& \frac{a_{N-2}}{6 dx_{N-1}} + \frac{a_{N-1} dx_{N-2}}{6 dx_{N-1}^{2}} , \end{matrix} \notag \\   \begin{matrix} 
   - \frac{2 a_{N-2} dx_{N-2}}{3 \left(dx_{N-2} + dx_{N-1}\right)^{2}} - \frac{a_{N-2}}{6 dx_{N-1}} - \frac{a_{N-1} dx_{N-2}}{6 dx_{N-1}^{2}}\end{matrix} \notag \\   \begin{matrix}  - \frac{a_{N} dx_{N-1}}{6 dx_{N}^{2}}, 
     \frac{a_{N} dx_{N-1}}{6 dx_{N}^{2}}\end{matrix}\bigg]  
\end{align}
The $N$ line
\begin{align}
  \bigg[\begin{matrix}\hdots ,& \frac{2 a_{N-1} dx_{N-1}}{3 \left(dx_{N-1} + dx_{N}\right)^{2}} & \frac{a_{N} dx_{N-1}}{6 dx_{N}^{2}} , \end{matrix} \notag \\   \begin{matrix}
     - \frac{2 a_{N-1} dx_{N-1}}{3 \left(dx_{N-1} + dx_{N}\right)^{2}} - \frac{a_{N} dx_{N-1}}{6 dx_{N}^{2}}\end{matrix} \bigg]
\end{align}

The reader will find an online environment for the numerical simulation at \url{https://colab.research.google.com/drive/1m6uhaur3eySqQ6eyjKf6SXxHXxXhSsWd?usp=sharing} and a git-hub depot of the numerical simulation at \url{https://github.com/christoautom/wave_1d}.
\end{extended}

\end{document}